%% file: main.tex
\documentclass[11pt,reqno]{amsart}

\usepackage{enumitem,kantlipsum}
\usepackage[table]{xcolor} 
\usepackage{booktabs} 
\usepackage{adjustbox} 
\usepackage{cancel}
\usepackage{epigraph}
\usepackage{geometry}                
\usepackage{subcaption}
\setlist[enumerate, 1]{label=\arabic*.  }

\usepackage{graphicx}
\usepackage{amssymb}
\usepackage{amsthm}
\usepackage{epstopdf}
\RequirePackage{doi}
\usepackage{hyperref}
\usepackage{graphicx}
\usepackage{forest}
\usepackage{tikz}
\usepackage{tikz-cd}
\usepackage{stmaryrd}
\usepackage{enumitem}
\usepackage{comment} 
\usepackage{newtxmath}
\usepackage{nicefrac}
\usepackage{pgfplots}
\pgfplotsset{compat=1.18}

\definecolor{USred}{cmyk}{0,1.00,0.65,0.34}
\renewcommand{\emph}[1]{{\textcolor{USred}{\em #1}}}

\hypersetup{
  colorlinks,
  citecolor=USred,
  linkcolor=USred,
  urlcolor=USred,
  allcolors=USred}

\def\SS{\mathbb{S}}

\def\GFR{\mathcal{R}} 

\DeclareMathOperator{\srec}{srec}
\DeclareMathOperator{\Flags}{RFlag}

\newcommand{\planted}{R_{\bullet}}

\theoremstyle{definition}
\newtheorem{de}{Definition}[section]
\theoremstyle{plain}
\newtheorem{thm}[de]{Theorem}

\newtheorem{lem}[de]{Lemma}
\newtheorem{p}[de]{Proposition}
\newtheorem{cor}[de]{Corollary}
\newtheorem{ex}[de]{Example}

\theoremstyle{definition}

\title{
On the enumeration of records of rooted trees and rooted forests.}
\author{ Adrián Lillo, Mercedes Rosas, Stefan Trandafir}

\begin{document}

\begin{abstract}
 A record of a rooted Cayley tree is a node whose label is the largest along the unique path to the root. In this work, we find elegant functional equations relating the generating functions for records of rooted Cayley trees and for records of forests of rooted trees with the Cayley tree function, and explore the consequences of our results. 
\end{abstract}
\maketitle

\input{tikz_preamble}

\vspace{-.2cm}

\input{1_Introduction}

\input{2_record_definitions}

\input{2B_record_decomposition}

\input{2C_enumerative_record_decomposition}

\input{3_generating_functions}

\input{4_record_numbers}
\input{5_consequences}
\input{6_final_comments}

\section*{acknowledgements}

This work has been partially supported by the Grant PID2020-117843GB-
I00 funded by MICIU/AEI/ 10.13039/501100011033. Author ST is supported by \href{10.3030/101070558}{FoQaCiA} which is funded by the European Union and NSERC.

\bibliographystyle{alpha}
\bibliography{references}

\end{document}

%% file: tikz_preamble.tex
\tikzset{
    record/.style={
        circle,
        draw=black,
        fill=black,
        inner sep=1.5pt,
        label={#1},
        solid
    },
    record/.default={},
    non-record/.style={
        circle,
        draw = black,
        fill = white,
        inner sep=1.5pt,
        solid
    },
    every label/.style={
            font=\tiny
        }
}

%% file: 1_Introduction.tex
\section{Introduction}

 A vertex of a rooted labelled tree $T$ is called a \emph{record} if its label is the largest along the path from the root to that vertex. The \emph{tree record numbers} count rooted  trees with a given number of vertices containing a specified number of records, while the \emph{forest record numbers} are defined analogously for rooted forests. The record numbers provide a refinement of Cayley’s classical formula for counting labelled trees by tracking the distribution of records. 
 We aim  to initiate a systematic study of the enumerative properties of the tree and forest record numbers. 
Despite their natural combinatorial significance, these numbers have received little attention in the existing literature.

 In their foundational book on Analytic Combinatorics \cite{Flajolet}, Flajolet and Sedgewick presented a summary of interesting counting results that can be expressed in terms of the \emph{Cayley tree function} $\mathcal{T}$. For example, they show that the exponential generating functions for mappings is 
 $\dfrac{1}{1 - \mathcal{T}}$, for connected mappings $\log \dfrac{1}{1 - \mathcal{T}}$, for mappings without fixed points $\dfrac{e^{-\mathcal{T}}}{1 - \mathcal{T}}$, for idempotent mappings  $e^{z e^z}$, for partial mappings $\dfrac{e^\mathcal{T}}{1 - \mathcal{T}}$, and so on.
Further interesting relations between the Cayley tree functions, its derivatives and related structures are described in \cite{Josual-V-tree-function} and
\cite{Cyril-tree-function}.

We present similarly elegant expressions for the \emph{tree record function}, the generating function for rooted trees according to their number of records, and for the \emph{forest record function}, the similarly defined generating function for rooted forests, expressed in terms of the Cayley tree function. Of course, these two expressions are related by the exponential formula.
\begin{align*}
&\text{Tree record function} & & 
\int_0^z  \frac{1}{s}
\frac{t \, \mathcal{T}(s)   }{1-  t \, \mathcal{T}(s)} \ ds &  \text{Proposition } \ref{le:trees_k_records}\\
&\text{Forest record function}  & & 
\exp \big(\int_0^z  \frac{1}{s}
\frac{t \, \mathcal{T}(s)   }{1-  t \, \mathcal{T}(s)} \ ds \big) \ \ \
& \text{Proposition } \ref{cor:R2B}
\end{align*}
Additional related expressions are  described in Section \ref{se:generating_functions}.

Our work rests on the \emph{record decomposition} of rooted labelled trees. In Section \ref{se:blob}, we provide an overview of this construction, as presented in \cite{LRT-Weary}. A closely related version of this construction has been described by Kreweras and Moszkowski \cite{Kreweras-Moszkowski}. 
In particular, we show in Lemma \ref{le:count_planted_fixed_type} that the number of rooted  trees with \emph{bonsai sequence} of type
$t$ is equal to
\begin{align*}
  \frac{1}{n}\, {n \choose t_1, t_2, \ldots, t_k} \,
  t_1^{t_1-1} t_2^{t_2-1} \ldots t_k^{t_k-1}.
\end{align*}
The precise definitions of a bonsai sequence and its type are given in Section \ref{se:blob}. Intuitively, the type of a bonsai sequence specifies the sizes of the trees that appear in its record decomposition.

In Proposition \ref{eq:forest_to_cayley_tree_functions} we describe some further relations between the record numbers of rooted trees and forests, and the Cayley tree function.
For example, we show that for  fixed $k$, the generating function for the numbers $R_\bullet(n, k)$ is given by 
\[
[t^k]\mathcal R_\bullet(z, t) = \frac{\mathcal{T}(z)^k}{k} - \frac{\mathcal{T}(z)^{k+1}}{k+1}. 
\]
 
The remainder of this article focuses directly on the record numbers. 
We show in Theorem \ref{prop:R_planted} that the \emph{tree record number}, counting rooted  trees with $n$ nodes and $k$ records, obeys the elegant closed formula:
    \begin{align}
    \label{Eq_rooted_introduction}
  k \,  (n-1)\ldots
    (n-k+1) \, n^{n-k-1}.
     \end{align}
Interestingly, as reported by the encyclopaedia of integer sequences (\url{https://oeis.org/A259334}), integers described by this formula have already appeared in the literature in the context of queues and traffic lights in the work of Haight \cite{Haight}. 
On the other hand, we also derive an explicit formula for the forest record numbers (Theorem \ref{thm:R(n,k) Stirling}). We show that the number of forests with $n$ nodes and $k$ records is given by
\[
\sum_{m = k+1}^{n+1} (-1)^{m+k-1} \binom{n}{m-1}m(n+1)^{n-m}c(m, m-k).
\]
where we denote by $c(\cdot,\cdot)$ the Stirling numbers of the first kind.

The rest of the paper deals with consequences of our results.
We describe two closely related recurrences, one for the tree record numbers and one for the forest record numbers (Proposition \ref{le:recurrence_forests}). They both rest on the record decomposition presented in Section \ref{se:blob}.
 We  describe the asymptotics of the tree record numbers (Corollary \ref{Cor:Lemma_k_asymptotics_records}).
We show that the exponential generating function for the tree record numbers (with $n \ge 1$ fixed) is given by the elegant polynomial expression $\frac{t}{n} (n+t)^{n-1}
$  (Corollary \ref{cor:record_fix_n}). 
We show that the sequences of  record numbers for rooted  trees obtained after fixing  $n$, and letting $k$ vary are log-concave (Corollary \ref{cor:log_concavity_trees}), and describe the location of the peaks (Corollary \ref{co:location_peaks_trees}).
On the other hand, we verify that the analogous sequences for the forest record numbers are log-concave up to $n=900$. 
Finally, in Table \ref{Fi:peak_fits}, we present data on the natural analogue of Kortchemski's statistic for permutations where we consider the depth of the record nodes.

  Our original motivation for studying records of rooted  trees and rooted forests comes from the following scenario.  Imagine some cars driving along a long one-way road, each with its own preferred speed. Overtaking is impossible, hence faster cars inevitably catch up to slower ones. If the slowest car happens to start first, it causes a major traffic jam. Even under better conditions,  cars tend to end up crawling behind slower vehicles.  
This clustering phenomenon is well known in queue theory (\cite{Haight}) and combinatorics (\cite{Foata}, \cite{Foata-Riordan}), where permutations represent car arrival orders,
cluster heads correspond to records of permutations (left-to-right maxima), and the unsigned Stirling numbers of the first kind count records in permutations.

In the sequel to this work \cite{LRT-bijective}, we present an elegant bijection, closely related to the one of Joyal \cite{Joyal}, which shows that the number of records in all rooted trees with $n$ nodes coincides with the number of connected endomorphisms of $[n]$. We also describe several sequences hidden in the table of tree record numbers, including well-known permutation statistics such as exceedances, small descents, and Kortchemski’s sums of the positions of the records. Furthermore, for the forest record numbers, we uncover connections to the weight of the symmetric group. Finally,  we ask: How many records does a random rooted  tree have? What is the probability of selecting a tree with $k$ records when we select uniformly at random a rooted  tree with $n$ nodes? What is the minimum $k$ for which a rooted  tree has at least $m$ nodes with probability at least $1/2$?, \&c. These questions have been considered in a different setting by Janson, \cite{Janson}. Our framework allows us to recover several of his results in a combinatorially explicit way.

The notion of records of rooted trees appears naturally in a closely related scenario in which cars travel along roads with merging junctions,  where the root indicates the observation point, and overtaking is not allowed. This scenario has been shown to be in bijection with parking functions in \cite{LRT-Weary}.  In addition, it was shown that the forest record numbers also count parking functions by the number of records, defined as the left-right maxima of the resulting parking configuration. 

Records of rooted  trees first appeared in work by Meir and Moon \cite{MM}, and more recently in the work of Janson \cite{Janson}, in relation to a pruning process on  trees. They have also been studied from a combinatorial point of view in \cite{MM, Kreweras-Moszkowski, Biane_J-V} and \cite{LRT-Weary}. However, despite the prominence of the study of trees in combinatorics, tree records remain widely underexplored. In this work, we attempt to examine them from an enumerative viewpoint.

%% file: 2_record_definitions.tex
\section{Records, record decomposition,  and record numbers.} \label{se:blob}

As all trees appearing in this work are assumed to be Cayley trees (finite labelled trees), we will simply refer to them as \emph{trees}. On the other hand, we refer to a collection of \emph{rooted trees} (with disjoint labels) as a \emph{rooted forest}. A node of a rooted tree is a \emph{record} if its label is the largest along the unique path to the root. In a rooted forest, a node is a \emph{record} if it is a record in its tree.

Let \emph{tree record number} $\planted(n,k)$ denote the number of rooted trees with $n$ nodes and $k$ records, and the \emph{forest record number} $R(n,k)$ be the number of  rooted forests with $n$ nodes and $k$ records.
It turns out to be convenient to consider a rooted forest $F$ as a single tree by introducing an additional node $\circ$, declaring it to be the root, and adding edges connecting it to the roots of all trees in $F$. 
In particular, if a rooted forest $F$ consists of a unique tree, this construction results in a planted tree, rooted at $\circ$.  
This construction not only unifies the treatment of records of rooted trees and of rooted forests in this work, but allows us to rely on results of \cite{LRT-Weary}. 

We emphasise that the term \emph{rooted tree} denotes a labelled tree, which can be rooted at any node, and that its root is always considered a record. By contrast, a \emph{tree rooted at \(\circ\)}, indicates that the tree \( T \) contains the distinguished node \(\circ\), which we assume to be its root but do not count as a record (nor as an actual node of the tree).
Therefore, the tree record numbers $\planted(n,k)$ also count the number of \emph{planted trees rooted at $\circ$}, with $n$ nodes other than $\circ$, and $k$ records, while the forest record numbers $R(n,k)$ also count the number of trees rooted at $\circ$, with $n$ nodes also excluding $\circ$,  and $k$ records. 

While $n>0$, both $R(n,0)$ and  $\planted(n, 0)$ are equal to zero,  as $n$ will always be a record, but $R(0,0) =1$, as the 
empty forest has no record, and $\planted(0, 0) = 0$, as the empty forest is not considered a tree. 

\begin{ex}
Figure \ref{fig:trees} illustrates all 16   trees with 3 vertices, each rooted at $\circ$. It also shows the 16 rooted forests on 3 nodes, where $\circ$ serves as a marker indicating the positions of the roots. We obtain the forest record numbers \( R(3,1) = 3 \),
 \( R(3,2) = 7 \), and 
\( R(3,3) = 6 \). On the other hand, restricting ourselves to the planted trees appearing in this figure, we obtain  the tree record numbers
$\planted(3,1) = 3$,
$\planted(3,2)= 4$, and
$\planted(3,3) = 2$.

 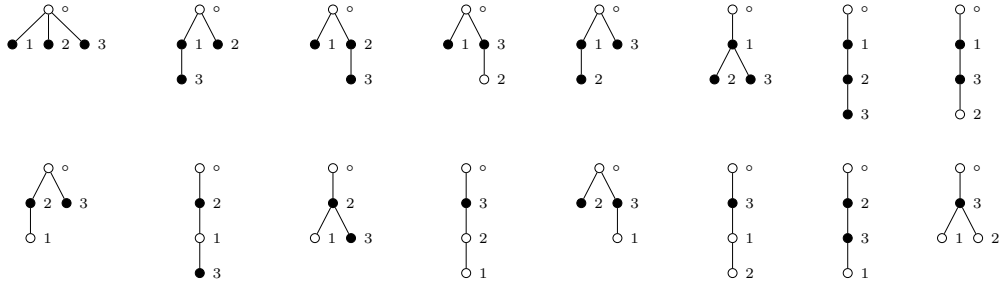
\begin{figure}[h!]
\centering
\scalebox{0.8}{
\normalfont
\input{all_trees_n4}
}
\caption{All 16  trees with 3 vertices  rooted at $\circ$. Record nodes are drawn in black. Nine of these trees are planted trees.}
\label{fig:trees}
\end{figure}
\end{ex}

We write  $[n]=\{1, 2, \ldots, n\}$ and $[n]_0=\{0,  1, 2, \ldots, n\}$. When we consider the extra root~$\circ$, we always label the root $\circ$ with 0.
 A rooted tree is said to be an \emph{increasing tree} if labels increase strictly along any path from the root to each of its leaves. Equivalently, if all its nodes are record nodes, with the exception of $\circ$, if it exists. Increasing trees have been widely studied in the combinatorial literature, see \cite{bergeron1992varieties, kuznetsov1994increasing}.

We present some values of the tree and forest records numbers Tables \ref{fig:planted_table} and \ref{fig:record_table}. The code used to generate this data  can be found at
\url{https://github.com/adrianlillo/tree_records}.
\input{table_tree_record_numbers}
\input{table_record_numbers}

%% file: all_trees_n4.tex
\usetikzlibrary{matrix}
\begin{tikzpicture}[
    level distance=0.5cm, 
    level 1/.style={sibling distance=0.6cm}, 
    level 2/.style={sibling distance=0.6cm},
    tree/.style={
        anchor=north, 
    }
]
\def\codeheight{15pt}
\matrix[matrix of nodes, row sep=0.5cm, column sep=1cm, nodes={tree}, ampersand replacement=\&] {
    \node[non-record, label={right:$\circ$}] (1-1-0) {}
    child {	
	    node[record, label={right:1}] (1-1-1) {}}
    child {	
	    node[record, label={right:2}] (1-1-2) {}}
    child {	
	    node[record, label={right:3}] (1-1-3) {}};\&
    \node[non-record, label={right:$\circ$}] (1-2-0) {}
    child {	
	    node[record, label={right:1}] (1-2-1) {}
	    child {	
		    node[record, label={right:3}] (1-2-3) {}}}
    child {	
	    node[record, label={right:2}] (1-2-2) {}};\&
    \node[non-record, label={right:$\circ$}] (1-3-0) {}
    child {	
	    node[record, label={right:1}] (1-3-1) {}}
    child {	
	    node[record, label={right:2}] (1-3-2) {}
	    child {	
		    node[record, label={right:3}] (1-3-3) {}}};\&
    \node[non-record, label={right:$\circ$}] (1-4-0) {}
    child {	
	    node[record, label={right:1}] (1-4-1) {}}
    child {	
	    node[record, label={right:3}] (1-4-3) {}
	    child {	
		    node[non-record, label={right:2}] (1-4-2) {}}};\&
    \node[non-record, label={right:$\circ$}] (2-1-0) {}
    child {	
	    node[record, label={right:1}] (2-1-1) {}
	    child {	
		    node[record, label={right:2}] (2-1-2) {}}}
    child {	
	    node[record, label={right:3}] (2-1-3) {}};\&
    \node[non-record, label={right:$\circ$}] (2-2-0) {}
    child {	
	    node[record, label={right:1}] (2-2-1) {}
	    child {	
		    node[record, label={right:2}] (2-2-2) {}}
	    child {	
		    node[record, label={right:3}] (2-2-3) {}}};\&
    \node[non-record, label={right:$\circ$}] (2-3-0) {}
    child {	
	    node[record, label={right:1}] (2-3-1) {}
	    child {	
		    node[record, label={right:2}] (2-3-2) {}
		    child {	
			    node[record, label={right:3}] (2-3-3) {}}}};\&
    \node[non-record, label={right:$\circ$}] (2-4-0) {}
    child {	
	    node[record, label={right:1}] (2-4-1) {}
	    child {	
		    node[record, label={right:3}] (2-4-3) {}
		    child {	
			    node[non-record, label={right:2}] (2-4-2) {}}}};\\ 
    \node[non-record, label={right:$\circ$}] (3-1-0) {}
    child {	
	    node[record, label={right:2}] (3-1-2) {}
	    child {	
		    node[non-record, label={right:1}] (3-1-1) {}}}
    child {	
	    node[record, label={right:3}] (3-1-3) {}};\&
    \node[non-record, label={right:$\circ$}] (3-2-0) {}
    child {	
	    node[record, label={right:2}] (3-2-2) {}
	    child {	
		    node[non-record, label={right:1}] (3-2-1) {}
		    child {	
			    node[record, label={right:3}] (3-2-3) {}}}};\&
    \node[non-record, label={right:$\circ$}] (3-3-0) {}
    child {	
	    node[record, label={right:2}] (3-3-2) {}
	    child {	
		    node[non-record, label={right:1}] (3-3-1) {}}
	    child {	
		    node[record, label={right:3}] (3-3-3) {}}};\&
    \node[non-record, label={right:$\circ$}] (3-4-0) {}
    child {	
	    node[record, label={right:3}] (3-4-3) {}
	    child {	
		    node[non-record, label={right:2}] (3-4-2) {}
		    child {	
			    node[non-record, label={right:1}] (3-4-1) {}}}};\&
    \node[non-record, label={right:$\circ$}] (4-1-0) {}
    child {	
	    node[record, label={right:2}] (4-1-2) {}}
    child {	
	    node[record, label={right:3}] (4-1-3) {}
	    child {	
		    node[non-record, label={right:1}] (4-1-1) {}}};\&
    \node[non-record, label={right:$\circ$}] (4-2-0) {}
    child {	
	    node[record, label={right:3}] (4-2-3) {}
	    child {	
		    node[non-record, label={right:1}] (4-2-1) {}
		    child {	
			    node[non-record, label={right:2}] (4-2-2) {}}}};\&
    \node[non-record, label={right:$\circ$}] (4-3-0) {}
    child {	
	    node[record, label={right:2}] (4-3-2) {}
	    child {	
		    node[record, label={right:3}] (4-3-3) {}
		    child {	
			    node[non-record, label={right:1}] (4-3-1) {}}}};\&
    \node[non-record, label={right:$\circ$}] (4-4-0) {}
    child {	
	    node[record, label={right:3}] (4-4-3) {}
	    child {	
		    node[non-record, label={right:1}] (4-4-1) {}}
	    child {	
		    node[non-record, label={right:2}] (4-4-2) {}
            }};\\};
\end{tikzpicture}

%% file: table_tree_record_numbers.tex
\begin{table}[h!]
  \renewcommand{\arraystretch}{1.3} 
\rowcolors{2}{gray!20}{white} 

    \begin{tabular}{c|ccccccccc}
        \toprule
        $n \backslash k$ & 1 & 2 & 3 & 4 & 5 & 6 & 7 & 8 & 9 \\
        \midrule
        1  & 1 & 0 & 0 & 0 & 0 & 0 & 0 & 0 & 0 \\
        2  & 1 & 1 & 0 & 0 & 0 & 0 & 0 & 0 & 0 \\
        3  & 3 & 4 & 2 & 0 & 0 & 0 & 0 & 0 & 0 \\
        4  & 16 & 24 & 18 & 6 & 0 & 0 & 0 & 0 & 0 \\
        5  & 125 & 200 & 180 & 96 & 24 & 0 & 0 & 0 & 0 \\
        6  & 1296 & 2160 & 2160 & 1440 & 600 & 120 & 0 & 0 & 0 \\
        7  & 16807 & 28812 & 30870 & 23520 & 12600 & 4320 & 720 & 0 & 0 \\
        8  & 262144 & 458752 & 516096 & 430080 & 268800 & 120960 & 35280 & 5040 & 0 \\
        \bottomrule
    \end{tabular}
        \caption{The tree record numbers $\planted(n,k).$} 
        \label{fig:planted_table}
\end{table}
\rowcolors{2}{}{} 

%% file: table_record_numbers.tex
\renewcommand{\arraystretch}{1.2} 
\rowcolors{2}{gray!20}{white} 

\begin{table}[h!]
  \renewcommand{\arraystretch}{1.3} 
\rowcolors{2}{gray!20}{white} 

    \begin{tabular}{c|cccccccccc}
        \toprule
        $n \backslash k$ & 1 & 2 & 3 & 4 & 5 & 6 & 7 & 8 
        \\
        \midrule
        1  & 1 & 0 & 0 & 0 & 0 & 0 & 0 & 0 
        \\
        2  & 1 & 2 & 0 & 0 & 0 & 0 & 0 & 0 
        \\
        3  & 3 & 7 & 6 & 0 & 0 & 0 & 0 & 0 
        \\
        4  & 16 & 39 & 46 & 24 & 0 & 0 & 0 & 0 
        \\
        5  & 125 & 310 & 415 & 326 & 120 & 0 & 0 & 0 
        \\
        6  & 1296 & 3240 & 4635 & 4360 & 2556 & 720 & 0 & 0 
        \\
        7  & 16807 & 42189 & 62825 & 65415 & 47656 & 22212 & 5040 & 0 
        \\
        8  & 262144 & 659456 & 1008448 & 1120385 & 927388 & 551852 & 212976 & 40320 
        \\
        \bottomrule
    \end{tabular}
        \caption{The forest record numbers $R(n,k).$} 
        \label{fig:record_table}
\end{table}
\rowcolors{2}{}{} 

%% file: 2B_record_decomposition.tex
\subsection{The record decomposition } \label{suse:record_decomposition}

The notion of a record a (labelled) tree naturally leads to the \emph{record decomposition} of \cite{LRT-Weary}. It is  closely related to the work of Kreweras and Moszkowski \cite{Kreweras-Moszkowski}.  We summarise it here.

\begin{figure}[hbt!] 
\centering
\input{tikz.record_cuts}
\caption{The tree $T$ with edges above record nodes deleted.}
\label{fi:record_deleted}
\end{figure}
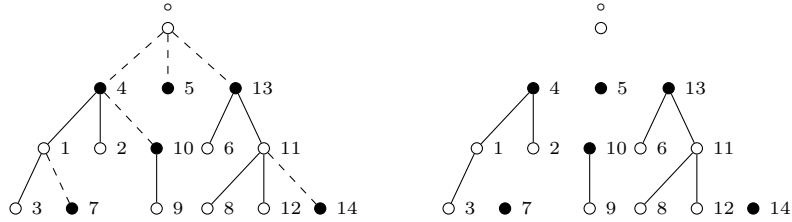

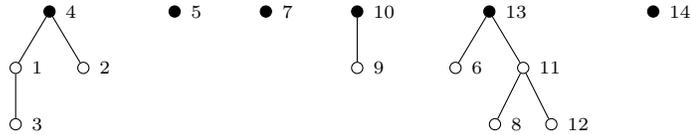
\begin{figure}[hbt!] 
\centering
\input{tikz.bonsai_sequence}
\caption{The record decomposition of $T$ has bonsai type $(4,1,1,2,5,1) \vDash 14$,  attachment sequence $(\circ, 1, 4, \circ, 11),$ and bonsai sequence as illustrated.}
\label{fi:bonsai_dec}
\end{figure}

A \emph{bonsai tree} is a rooted   tree that has its root as its unique record. Observe that this this implies that a   tree rooted at $\circ$ is never a bonsai because $\circ$ is never a record. 
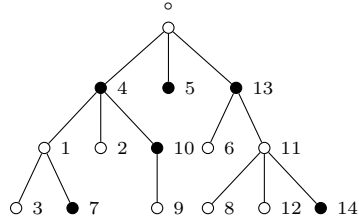
\begin{figure}[hbt!] 
\centering
\input{example_records}
\caption{A tree $T$ labelled with \( [14]_0 \) and rooted at $\circ$. The record nodes of $T$ are drawn in black. The root \( \circ \) is not considered a record node. 
}\label{Running_Example}
\end{figure}
Let \( T \) be a   tree rooted at  \( \circ \). The \emph{bonsai sequence $b(T)$} of \( T \) is the sequence of bonsai trees \( b(T) = (B_1, B_2, \ldots, B_k) \) obtained by deleting from \( T \) all edges that connect a record node to its parent (this operation always yields a rooted forest of bonsai trees; see Figure~\ref{fi:record_deleted}), and ordering the resulting bonsais increasingly according to the labels of their roots, as illustrated in Figure~\ref{fi:bonsai_dec}.

A \emph{composition} of $n$, writen \emph{$t \vDash n$}, is a way of writing $n$ as the sum of a sequence of \emph{positive integers}. The number of integers appearing in the sequence is called its \emph{length} and denoted by \emph{$\ell(t)$}. 
The \emph{bonsai type} $ t(T)$ of $T$ is defined as \( t(T) = (|B_1|, |B_2|, \ldots, |B_k|) \). Observe that  if \( t(T) \vDash n \), then the number of nodes in \( T \) excluding the root is equal to $n$, and  $\ell(t(T))$ is the number of records in $T$.

Different trees can share the same bonsai sequence $b = (B_1, B_2, \ldots, B_k)$. If we tried to reconstruct a tree from $b$, we would observe that, while bonsai $B_1$ must be attached to $\circ$, bonsai $B_2$ could be attached to any of the nodes of $B_1$ or to $\circ$, bonsai $B_3$ could be attached to any node of $B_1$, $B_2$ or to $\circ$, and so on.  

However, if in addition to the bonsai sequence of $T$, we record the parent of each record node,  we can fully reconstruct the tree $T$ from $b$. Moreover, as the parent of the first node is redundant (it is always attached to $\circ$) we can omit this information. The sequence of parents of record nodes  (other than the first one) is called the \emph{attachment sequence} of $T$. 
The attachment sequence of our running example appears in Figure \ref{fi:bonsai_dec}. 

The \emph{record decomposition} of a   tree $T$ rooted at $\circ$ is the ordered pair consisting of its bonsai sequence and its attachment sequence.
The enumerative results presented in this work rely heavily on the record decomposition, particularly on the following three lemmas.

%% file: example_records.tex
\normalfont
\begin{tikzpicture}[
    level distance=0.8cm,
    level 1/.style={sibling distance=0.9cm},
    level 2/.style={sibling distance=0.75cm,},
    level 3/.style={sibling distance=0.75cm},
    ]
    \node[non-record, label={$\circ$}] (1) {}
    child {
    	node[record, label={right:4}] (5) {}
    	child {
    		node[non-record, label={right:1}] (2) {}
    		child {
    			node[non-record, label={right:3}] (4) {}
    		}
    		child {
    			node[record, label={right:7}] (8) {}
    		}
    	}
    	child {
    		node[non-record, label={right:2}] (3) {}
    	}
    	child {
    		node[record, label={right:10}] (11) {}
    		child {
    			node[non-record, label={right:9}] (10) {}
    		}
    	}
    }
    child {
    	node[record, label={right:5}] (6) {}
    }
    child {
    	node[record, label={right:13}] (14) {}
    	child {
    		node[non-record, label={right:6}] (7) {}
    	}
    	child {
    		node[non-record, label={right:11}] (12) {}
    		child {
    			node[non-record, label={right:8}] (9) {}
    		}
    		child {
    			node[non-record, label={right:12}] (13) {}
    		}
    		child {
    			node[record, label={right:14}] (15) {}
    		}
    	}
    };
\end{tikzpicture}

%% file: tikz.record_cuts.tex
\begin{tikzpicture}
    \matrix(outer) [matrix of nodes, row sep=0.3cm, column sep=1cm, ampersand replacement=\&]{ 
    \scoped
    [
    level distance=0.8cm,
    level 1/.style={sibling distance=0.9cm},
    level 2/.style={sibling distance=0.75cm,},
    level 3/.style={sibling distance=0.75cm},
    cutedge/.style={dashed},
    ]
    \node[non-record, label={$\circ$}] (0) {}
    child {node[record, label={right:4}] (4) {}
        child {node[non-record, label={right:1}] (1) {}
            child {node[non-record, label={right:3}] (3) {}
                edge from parent[solid]}
            child {node[record, label={right:7}] (7) {}
                edge from parent[cutedge]}
            edge from parent[solid]}
        child {node[non-record, label={right:2}] (2) {}
                edge from parent[solid]}
        child {node[record, label={right:10}] (10) {}
            child {node[non-record, label={right:9}] (9) {}
                    edge from parent[solid]}
            edge from parent[cutedge]}
        edge from parent[cutedge]
    }
    child {node[record, label={right:5}] (5) {}
        edge from parent[cutedge]}
    child {node[record, label={right:13}] (13) {}
        child {node[non-record, label={right:6}] (6) {}
                edge from parent[solid]}
        child {node[non-record, label={right:11}] (11) {}
            child {node[non-record, label={right:8}] (8) {}
                    edge from parent[solid]}
            child {node[non-record, label={right:12}] (12) {}
                    edge from parent[solid]}
            child {node[record, label={right:14}] (14) {}
                    edge from parent[cutedge]}
                edge from parent[solid]
        }
            edge from parent[cutedge]
        }; \&
    \scoped 
    [
    level distance=0.8cm,
    level 1/.style={sibling distance=0.9cm},
    level 2/.style={sibling distance=0.75cm,},
    level 3/.style={sibling distance=0.75cm},
    cutedge/.style={draw=none},
    ]
    \node[non-record, label={$\circ$}] (0) {}
    child {node[record, label={right:4}] (4) {}
        child {node[non-record, label={right:1}] (1) {}
            child {node[non-record, label={right:3}] (3) {}
                edge from parent[solid]}
            child {node[record, label={right:7}] (7) {}
                edge from parent[cutedge]}
            edge from parent[solid]}
        child {node[non-record, label={right:2}] (2) {}
                edge from parent[solid]}
        child {node[record, label={right:10}] (10) {}
            child {node[non-record, label={right:9}] (9) {}
                    edge from parent[solid]}
            edge from parent[cutedge]}
        edge from parent[cutedge]
    }
    child {node[record, label={right:5}] (5) {}
        edge from parent[cutedge]}
    child {node[record, label={right:13}] (13) {}
        child {node[non-record, label={right:6}] (6) {}
                edge from parent[solid]}
        child {node[non-record, label={right:11}] (11) {}
            child {node[non-record, label={right:8}] (8) {}
                    edge from parent[solid]}
            child {node[non-record, label={right:12}] (12) {}
                    edge from parent[solid]}
            child {node[record, label={right:14}] (14) {}
                    edge from parent[cutedge]}
                edge from parent[solid]
        }
            edge from parent[cutedge]
        }; \\
        };
\end{tikzpicture}

%% file: tikz.bonsai_sequence.tex
\begin{tikzpicture}[
    level distance=0.75cm,
    level 1/.style={sibling distance=0.9cm},
    level 2/.style={sibling distance=0.75cm},
    level 3/.style={sibling distance=0.75cm},
]
\matrix[matrix of nodes, row sep=0.3cm, column sep=0.6cm, ampersand replacement=\&]
{
    \node[record, label={right:4}] (4) {}
        child {node[non-record, label={right:1}] (1) {}
            child {node[non-record, label={right:3}] (3) {}
                edge from parent[solid]}
            edge from parent[solid]}
        child {node[non-record, label={right:2}] (2) {}
                edge from parent[solid]}; \&
    \node[record, label={right:5}] (5) {};\&
    \node[record, label={right:7}] (7) {};\&
    \node[record, label={right:10}] (10) {}
        child {node[non-record, label={right:9}] (9) {}
                edge from parent[solid]};\&
    \node[record, label={right:13}] (13) {}
        child {node[non-record, label={right:6}] (6) {}
                edge from parent[solid]}
        child {node[non-record, label={right:11}] (11) {}
            child {node[non-record, label={right:8}] (8) {}
                    edge from parent[solid]}
            child {node[non-record, label={right:12}] (12) {}
                    edge from parent[solid]}
        };\&
    \node[record, label={right:14}] (14) {};\\
};
\end{tikzpicture}

%% file: 2C_enumerative_record_decomposition.tex
\subsection{Some enumerative consequences of the record decomposition.}
We count the number of rooted  trees and rooted forests, according to the sizes of the bonsais appearing in its record decomposition. 
As described in the introduction, we  model rooted trees and rooted forests by  rooted trees rooted at $\circ$.

A \emph{restricted flag $\mathcal{F}$ of $[n]$} of length $k$  is a sequence of subsets
$
    S_1 \subseteq S_2 
    \subseteq \ldots \subseteq S_{k-1}
    \subseteq S_k=[n]
$
such that the maximum element of $S_{i+1}$  does not belong to $S_i$, for $i= 1, \ldots, k-1$. We set $S_0=\emptyset$ and $S_k=[n]$. 
Observe that any bonsai sequence $b=(B_1, B_2, \ldots, B_k)$ defines a restricted flag: set $S_i$ to be the set of labels of $\bigcup_{j=1}^{i} B_j$, for $ 1\le i \le n.$ 

A restricted flag $\mathcal{F}$ has \emph{type} 
$t=(t_1, t_2, \ldots, t_k)$
if  for each $i\ge 1$, $t_i = \big| S_{i+1} \setminus S_{i} \big|$.
Note that $t\vDash n$ is a composition of $n$ of length $k$, and that no part of $t$ is equal to zero.

\begin{ex}
The restricted flag of the bonsai sequence of our running example \ref{fi:bonsai_dec} is:
\begin{multline*}
\emptyset \subseteq
\{1, 2, 3, 4 \} 
\subseteq
\{1, 2, 3, 4, 5 \}
\subseteq
\{1, 2, 3, 4, 5, 7 \}\\
\subseteq
\{1, 2, 3, 4, 5, 7, 9, 10 \}
\subseteq
\{1, 2, 3, 4, 5, 6,  7, 8, 9, 10, 11, 12, 13 \}\\
\subseteq
\{1, 2, 3, 4, 5, 6,  7, 8, 9, 10, 11, 12, 13,14 \}.
\end{multline*}
Its  type is 
 $(4, 1,1, 2, 5, 1)\vDash 14.$ Observe that the differences between successive sets record their labels.
\end{ex}
Given a non-negative integer $i$, and a composition $t=\{t_1, t_2, \cdots, t_n\}$,
define the symbols  
\emph{$[t]_i= t_1 + t_2+\cdots+t_i$},  and
\emph{$[t]_i!= [t]_1 [t]_2 \ldots [t]_i$}.
This notation is close to the one used for the $q$-integer and $q$ factorial coefficients, which can be easily recovered from them by specialising in $t_i=q^{i-1}$.

\begin{lem}\label{le:restricted_flags}
Let $t=(t_1, t_2, \ldots, t_k)$ be a composition of $ n$.
The number of restricted flags of type 
$t$ is given by
\begin{align*}
\Flags(t)  
= 
{n \choose t_1, t_2, \ldots, t_k} \frac{t_1 t_2 \ldots t_k}{[t]_k!}
= \frac{1}{n}
{n \choose t_1, t_2, \ldots, t_k} \frac{t_1 t_2 \ldots t_k}{[t]_{k-1}!}
\end{align*}

\end{lem}

Note that the number of restricted flags of type $t$ is {\bf not} invariant under permutations of the parts of $t$.

\begin{proof}
We construct the restricted flag \( \mathcal{F} \) from top to bottom, selecting  at step $i$ a subset \( S_{k-i} \) of $S_{k-i+1}$, of cardinality $t_1+t_2+\ldots+t_{k-i}$ and with the constraint that we are not allowed to choose the largest element in \( S_{k-i+1} \). Thus,
\[
\Flags(t) = 
{t_1 + t_2 + \ldots +t_{k} -1
\choose
t_1 + t_2 + \ldots +t_{k-1}
}
{t_1 + t_2 + \ldots +t_{k-1} -1
\choose
t_1 + t_2 + \ldots + t_{k-2}
}
\ldots 
{t_1+t_2 + t_3-1 \choose t_1+t_2}
{t_1+t_2-1 \choose t_1}.
\]
Using the definition of $[t]_k$, and the symmetry of the binomial coefficients, this identity can be rewritten as
\[
\Flags(t) = 
{[t]_k -1
\choose
t_k-1
}
{[t]_{k-1} -1
\choose
t_{k-1}-1
}
\ldots 
{[t]_3 -1
\choose
t_3-1
}
{[t]_2 -1
\choose
t_2-1
}
\]
Therefore,
 applying the identity 
$
\frac{k}{n} {n \choose k} = {n-1 \choose k-1}
$ to each binomial coefficient in the expression, and using the usual expression of the multinomial coefficients as a product of binomial coefficients, we obtain that
\begin{align*}
\Flags(t)  
=
\frac{t_k}{[t]_k}
\frac{t_{k-1}}{[t]_{k-1}}
\cdots
\frac{t_3}{[t]_3}
\frac{t_2}{[t]_2}
\frac{t_1}{[t]_1}
{[t]_k \choose t_1, t_2, \ldots, t_k}
=
{n \choose t_1, t_2, \ldots, t_k} \frac{t_1 t_2 \ldots t_k}{[t]_{k}!}.
\end{align*}
where we recall that $n = [t]_k = t_1+t_2+\ldots+t_{k}$, and observe that $[t]_1/t_1=1$.

\end{proof}

Now, we proceed to count the number of bonsai sequences with $\mathcal{F}$ as its restricted flag.  For each $i$, we construct a bonsai tree with $t_i$
nodes, with labels in $S_i\setminus S_i-1$, and rooted at its largest element, $r_i$. 
Cayley's formula tells us that the number of ways in which this can be achieved is
$
t_i^{t_i-2}
$. 
We iterate this construction for $i=1, \ldots, k$, obtaining a sequence $b=(B_1, B_2, \ldots, B_k)$ of bonsais. 
Moreover, from  the definition of restricted flag, $r_1 \le r_2 \le \ldots \le r_k$ and $b$ is a bonsai sequence. We record this observation in Lemma \ref{le:bonsais_flag}.

\begin{lem} \label{le:bonsais_flag}
The number of bonsai sequences that can be constructed from the restricted flag $\mathcal{F}$ of type 
$t=(t_1, t_2, \ldots, t_k) \vDash n$ is
\[
 t_1^{t_1-2} t_2^{t_2-2} \ldots t_k^{t_k-2}
 \]
 \end{lem}

Finally, we take into consideration the potential number of attachment sequences. 
 
\begin{lem} \label{lem:same_Bonsai_seq}
Let $b=(B_1, B_2, \ldots, B_k)$ be a bonsai sequence of type $(t_1, t_2, \ldots, t_k)$. Then, the number of rooted forests with bonsai sequence $b$ is 
\begin{align}\label{eq:same_Bonsai_seq}
(t_1 + 1) (t_1+t_2 + 1) \ldots (t_1+t_2 + \ldots + t_{k-1} + 1)  
\end{align}
\end{lem}

\begin{proof} Recall that we model a rooted forest 
by a tree rooted at $\circ$.
Then we have to attach bonsai $B_1$ to $\circ,$ bonsai $B_2$ to any element of $S_1 \cup \circ$, bonsai $B_3$ to any element of $S_2 \cup \circ$, etc. 

 \end{proof}

  Remarkably, applying the specialization $x_i\mapsto q^i$ to Equation (\ref{eq:same_Bonsai_seq}), we obtain the $q$-factorial coefficients, the generating function for permutations of $\mathbb S_n$ by the number of inversions. 

A small variation of this argument implies the following result for planted trees/rooted trees.

\begin{lem} \label{lem:same_Bonsai_seq_planted}
Let $b=(B_1, B_2, \ldots, B_k)$ be a bonsai sequence of type $(t_1, t_2, \ldots, t_k)$. Then, the number of rooted  trees with bonsai sequence $b$ is 
\begin{align}\label{eq:same_Bonsai_seq_planted}
t_1 (t_1+t_2 ) \cdots (t_1+t_2 + \ldots + t_{k-1} )
=[t]_{k-1}!
\end{align}
\end{lem}
\begin{proof}
Recall that we model a rooted  tree by a planted tree rooted at $\circ$.
    Since we cannot graft any bonsai (other than the first) at $\circ$, the number of choices for each entry of the attachment sequence decreases by one with respect to Lemma \ref{lem:same_Bonsai_seq}.
\end{proof}

\begin{lem}
\label{le:count_planted_fixed_type}
Let $t = (t_1, t_2, \ldots, t_k) \vDash n$.
The number of rooted trees labelled with $[n]$ and with bonsai type
$t$ is equal to
\begin{align}
\label{tree_record_composition}
  \frac{1}{n}\, {n \choose t_1, t_2, \ldots, t_k} \,
  t_1^{t_1-1} t_2^{t_2-1} \cdots t_k^{t_k-1}.
\end{align}
Note that this expression does not depend on the order of the parts of $t$, but only on the underlying partition. 
\end{lem}

Equation (\ref{tree_record_composition}) also counts the number of planted trees labelled with $[n]_0$, rooted at $\circ$, and with bonsai type $t$.

\begin{proof}
The record decomposition implies that rooted trees are in bijection with triples consisting of a restricted flag of type $t$ (counted in Lemma \ref{le:restricted_flags}), a set of trees with that precise set of labels (counted in Lemma \ref{le:bonsais_flag}), and an attachment sequence
(counted in Lemma \ref{lem:same_Bonsai_seq_planted}). Multiplying the expressions obtained in these three lemmas, we obtain 
\begin{align}
\label{eq_no_simpl}
\binom{n - 1}{t_1 - 1, t_2, \dots,t_k} \,
  t_1^{t_1-2} t_2^{t_2-1} \cdots t_k^{t_k-1} =
  \frac{1}{n}\, {n \choose t_1, t_2, \ldots, t_k} 
  \,
  t_1^{t_1-1} t_2^{t_2-1} \cdots t_k^{t_k-1}. 
\end{align}\qedhere
\end{proof}
While it is tempting to sum over all compositions of $n$ with $k$ parts to obtain a formula for the tree record numbers, a more elegant and useful result is presented in Theorem~\ref{prop:R_planted}.
On the other hand,  a similar argument, this time based on Lemmas \ref{le:restricted_flags},  \ref{le:bonsais_flag}, and \ref{lem:same_Bonsai_seq}, yields a comparable but less elegant formula for the enumeration of rooted forests.

\begin{p}
\label{thm:record_count}
\label{le:count_forests_fixed_type}
Let $t = (t_1, t_2, \ldots, t_k) \vDash n$.
The number of  rooted forests of bonsai type $t$ is equal to
\begin{align}
\label{forest_record_composition}
 \frac{1}{n }  \,
 {n \choose t_1, t_2, \ldots, t_k}
   \,
  t_1^{t_1-1} t_2^{t_2-1} \cdots t_k^{t_k-1}
  \prod_{i=i}^{k-1}
  \left(1+ \frac{1}{[t]_i}\right)
\end{align}
Moreover,
the  number of rooted forests with $k$ records is equal to
\begin{align*}
\label{le:count_forests}
R(n,k) = \frac{1}{n }  \,
\sum_{\substack{t \vDash n \\ \ell(t) = k}}
 {n \choose t_1, t_2, \ldots, t_k}
  t_1^{t_1-1} t_2^{t_2-1} \cdots t_k^{t_k-1}
  \prod_{i=i}^{k-1}
  \left(1+ \frac{1}{[t]_i}\right)
\end{align*}

\end{p}

Equation (\ref{forest_record_composition}) also counts the number of trees labelled with $[n]_\circ$, rooted at $\circ$, and  with bonsai type $t$.

Finally, using a partial version of the record decomposition, we obtain recursive formulas for both families of record numbers. 
 
\begin{p}
\label{le:recurrence_trees}
\label{le:recurrence_forests}
The tree record number $\planted(n, k)$ can be defined recursively by
\begin{align*}
\planted(n, k) &= \sum_{i = k-1}^{n-1}\binom{n-1}{i}
\planted(i, k-1)(n-i)^{n-i-2}\ i, \qquad n\geq 1, \qquad 2 \leq k \leq n;
\\
\planted(n, 1) &= n^{n-2},  \quad n \geq 1; \qquad \planted(n, k) = 0 \quad \text{otherwise}.
\end{align*} 
On the other hand, the forest record numbers $R(n,k)$ satisfy the recursion \begin{align*}
R(n, k) &= \sum_{i = k-1}^{n-1}\binom{n-1}{i}R(i, k-1)(n-i)^{n-i-2}(i+1), \qquad n\geq1, \qquad 1\leq k \leq n; \\
R(0, 0) &= 1; \qquad R(n,k) = 0 \quad \text{otherwise}.
\end{align*}
\end{p}

\begin{figure}
    \centering
    \resizebox{0.75\textwidth}{!}{$%
    \vcenter{\hbox{\input{tikz.record_recursion_A}}}%
  \hspace{0.5cm}%
  \equiv%
  \hspace{0.85cm}%
  \vcenter{\hbox{\input{tikz.record_recursion_B}}}%
$}
    \caption{The recursive decomposition of Lemma \ref{le:recurrence_forests}.}
    \label{fig:recursion}
\end{figure}
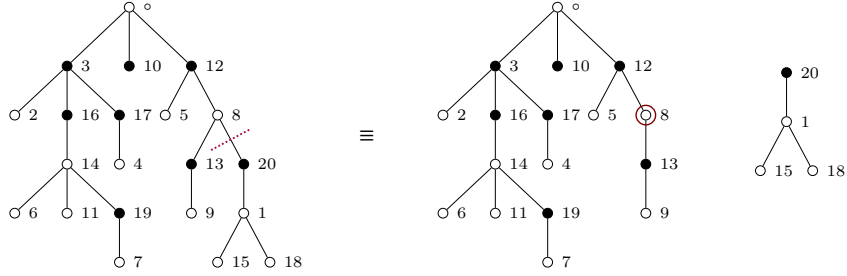

\begin{proof} To be able to prove both results simultaneously, we let $T$ be a  tree rooted at $\circ$, labelled with $\circ, 1, 2, \ldots, n$, and with $k$ records, and recall that for the tree record numbers we have to consider the additional condition that $T$ is a planted tree.

Delete from $T$ the edge joining $n$ with its parent $p$. This results in three pieces of information (see Figure \ref{fig:recursion}):
\begin{enumerate}[wide]
     
\item[i)] A tree $T'$ rooted at $\circ$, with $i$ non-root nodes and $k-1$ records. 
There are ${n-1 \choose i}$ possibilities for the set of labels $S$. On the other hand, there are $\planted(i,k-1)$ rooted trees (or $R(i, k-1)$ rooted forests) that we root at $\circ$, label with $S$, and that have $k-1$ records.

\item[ii)] A bonsai rooted at $n$, with $n-i$ nodes. The bonsai tree is labelled with the entries in $\bar S$, the complement of $S$ in $[n]$. There are $(n-i)^{n-i-2}$ unrooted trees in this set of labels, that we always root at $n.$

\item[iii)] A distinguished vertex $p$ of $T'$. 
There are two different situations to consider here.
In the case of rooted trees, $p$ can be any node of $T'$ other that $\circ$. Hence, it can be chosen in $i$ different ways. On the other hand, for rooted forests, $p$ can be any node of $T'$, including $\circ$. Hence, it can be chosen in $i+1$ different ways.\vspace{0.1cm}
\end{enumerate}
 Multiplying these numbers and summing over all possible $i$ we obtain our desired result.
\end{proof}

%% file: tikz.record_recursion_A.tex
\begin{tikzpicture}
[
level 1/.style={sibling distance = 0.95cm, level distance=0.9cm},
level 2/.style={sibling distance = 0.8cm, level distance=0.75cm},
sibling distance = 0.9cm,
]
\node[non-record, label={right:$\circ$}] (N-0) {}
child {	node[record, label={right:3}] (N-3) {}
	child {	node[non-record, label={right:2}] (N-2) {}}
	child {	node[record, label={right:16}] (N-16) {}
		child {	node[non-record, label={right:14}] (N-14) {}
			child {	node[non-record, label={right:6}] (N-6) {}}
			child {	node[non-record, label={right:11}] (N-11) {}}
			child {	node[record, label={right:19}] (N-19) {}
				child {	node[non-record, label={right:7}] (N-7) {}}}}}
	child {	node[record, label={right:17}] (N-17) {}
		child {	node[non-record, label={right:4}] (N-4) {}}}}
child {	node[record, label={right:10}] (N-10) {}}
child {	node[record, label={right:12}] (N-12) {}
	child {	node[non-record, label={right:5}] (N-5) {}}
	child {	node[non-record, label={right:8}] (N-8) {}
		child {	node[record, label={right:13}] (N-13) {}
			child {	node[non-record, label={right:9}] (N-9) {}}}
		child {	node[record, label={right:20}] (N-20) {}
			child {	node[non-record, label={right:1}] (N-1) {}
				child {	node[non-record, label={right:15}] (N-15) {}}
				child {	node[non-record, label={right:18}] (N-18) {}}}}}};
\newcommand{\cutlen}{0.8}

\coordinate (M) at ($(N-20)!0.5!(N-8)$);   
\coordinate (P1) at ($(M)!\cutlen!90:(N-8)$); 
\coordinate (P2) at ($(M)!\cutlen!-90:(N-8)$);

\draw[thick, USred, densely dotted] (P1) -- (P2);
\end{tikzpicture}

%% file: tikz.record_recursion_B.tex
\begin{tikzpicture}
[
level 1/.style={sibling distance = 0.95cm, level distance=0.9cm},
level 2/.style={sibling distance = 0.8cm, level distance=0.75cm},
sibling distance = 0.8cm,
]

\node[non-record, label={right:$\circ$}] (N-0) {}
child {	node[record, label={right:3}] (N-3) {}
	child {	node[non-record, label={right:2}] (N-2) {}}
	child {	node[record, label={right:16}] (N-16) {}
		child {	node[non-record, label={right:14}] (N-14) {}
			child {	node[non-record, label={right:6}] (N-6) {}}
			child {	node[non-record, label={right:11}] (N-11) {}}
			child {	node[record, label={right:19}] (N-19) {}
				child {	node[non-record, label={right:7}] (N-7) {}}}}}
	child {	node[record, label={right:17}] (N-17) {}
		child {	node[non-record, label={right:4}] (N-4) {}}}}
child {	node[record, label={right:10}] (N-10) {}}
child {	node[record, label={right:12}] (N-12) {}
	child {	node[non-record, label={right:5}] (N-5) {}}
	child {	node[non-record, label={right:8}] (N-8) {}
		child {	node[record, label={right:13}] (N-13) {}
			child {	node[non-record, label={right:9}] (N-9) {}}}
		}};

\node[semithick, circle, draw, red!50!black, inner sep=3pt] (c) at (N-8){};

\begin{scope}[level 1/.style={sibling distance = 0.8cm, level distance=0.75cm}]
\node[record, label={right:20}] (M-9) at (3.5,-1) {}
child {	node[non-record, label={right:1}] (M-2) {}
	child {	node[non-record, label={right:15}] (M-1) {}}
	child {	node[non-record, label={right:18}] (M-5) {}}};
\end{scope}
\end{tikzpicture}

%% file: 3_generating_functions.tex
\section{The Cayley tree function, and record generating functions for rooted trees and rooted forests}
\label{se:generating_functions}

In this section, we study the relationship between the generating function for records of trees, records of forests, and the ubiquitous Cayley tree function. 
\begin{align*}
& \mathcal{T}(z) 
=  \sum_{n\ge1} n^{n-1}\, \frac{z^{n} }{n!},
&&\text{\emph{Cayley tree function.}}\\
&\GFR_\bullet(z, t) = \sum_{n,k\ge 0}  \planted(n,k) \frac{z^{n}}{n!} t^k, 
&&\text{\emph{tree record function}}\\
&\GFR(z, t) = \sum_{n,k\ge 0}  R(n,k) \frac{z^{n}}{n!} t^k, 
&&\text{\emph{forest record function}}
\end{align*}
Let $A$ be the ring of polynomials with real coefficients.We consider our generating functions as elements of $A[[z]]$. 
Indeed, the tree and the forest record generating functions begin as follows:
\[
\begin{split}
\mathcal R_\bullet(z,t) = t z &+ \left(t +  t^{2}\right) \frac{z^{2}}{2!} + \left(3 t + 4 t^{2} + 2 t^{3}\right) \frac{z^{3}}{3!} \\ &+ \left(16 t + 24 t^{2} + 18 t^{3} + 6 t^{4}\right) \frac{z^{4}}{4!} + \left(125 t + 200 t^{2} + 180 t^{3} + 96 t^{4} + 24 t^{5}\right) \frac{z^{5}}{5!} + \cdots
\end{split}
\]
\[
\begin{split}
\mathcal R(z, t) = 1 &+ t z + \left(t + 2 t^{2}\right) \frac{z^{2}}{2!} + \left(3 t + 7 t^{2} + 6 t^{3}\right) \frac{z^{3}}{3!} \\ &+ \left(16 t + 39 t^{2} + 46 t^{3} + 24 t^{4}\right) \frac{z^{4}}{4!} + \left(125 t + 310 t^{2} + 415 t^{3} + 326 t^{4} + 120 t^{5}\right) \frac{z^{5}}{5!} + \cdots
\end{split}
\]

Before starting our analysis of our record generating functions, we recall some key properties of the Cayley tree function. First, it can be defined implicitly by the functional equation:
\[
\mathcal{T}(z) = z \exp\left( \mathcal{T}(z) \right).
\]
This classical result follows directly from the recursive structure of rooted trees: removing the root of a tree yields a set of rooted trees. For a detailed discussion of the Cayley tree function, see page~127 of~\cite{Flajolet}.
Alternatively, the Cayley tree function can be expressed directly in terms of the Lambert~$W$-function as 
$
\mathcal{T}(z) = -W(-z),
$ see ~\cite{Flajolet}.

The exponential formula gives a functional equation relating the tree record function with the forest record function.  After all, a rooted forest is a set of rooted trees, and the number of records is additive. 
\begin{lem}
\label{le:Exponential_formula_record_gf}
The tree record  function and the  forest record function are related by the functional equation:
\begin{align}
\mathcal{R}(z, t) = \exp({\mathcal{R}_\bullet(z,t))}.
\end{align}
\end{lem}
Our next results describe a functional equation that connects the Cayley tree function and the tree record function. Before presenting it, we provide two illustrative examples.
\begin{ex}
\label{ex_ge_unrooted} We start with a familiar example.
The generating function for \emph{unrooted  trees} with at least one node is given by the integral
\begin{align*}
\int_0^z  \frac{1}{s}
   \mathcal{T}(s)  \ ds.
   \end{align*} Indeed, this is a straightforward calculation:
   \begin{align*}
\int_0^z  \frac{1}{s}
   \left( \sum_{n\ge1} n^{n-2} \frac{s^n}{(n-1)!} \right) \ ds 
  &=  
   \sum_{n\ge1} n^{n-2} \int_0^z  \frac{s^{n-1}}{(n-1)!}  \ ds 
   =  
   \sum_{n\ge1} n^{n-2}  \frac{z^{n}}{ n!}  
\end{align*}
Since there are $n^{n-2}$ unrooted tree on $n$ vertices, we have obtained the generating function for unrooted  trees.
In particular,  this result implies that the Cayley tree function satisfies the relation
\[
\int_0^z  \frac{
    \mathcal{T}(s)}{s}  \ ds
   = T(z) - T^2(z).
\]
\end{ex}

\begin{ex}
The generating function for \emph{rooted trees with two records} is given by
\begin{align*}
\int_0^z  \frac{1}{s}
  \big(t \, \mathcal{T}(s) \big)^2 \ ds
  \end{align*} 
  Let us verify this claim. The previous expression is equal to
  \allowdisplaybreaks[1]
  \begin{align*}
\int_0^z  \frac{1}{s}
   \left(\sum_{n\ge1} t \ n^{n-1} \frac{s^n}{n!}\right)^2  \ ds 
  &=  
  t^2 \sum_{n\ge1}\sum_{\substack{n_1, n_2 \ge1\\n_1+n_2=n}} n! \ \frac{n_1^{n_1-1}}{n_1!}\frac{n_2^{n_2-1}}{n_2!}   \int_0^z  \frac{s^{n-1}}{ n!}  \ ds \\
   &=  
  t^2 
  \sum_{n\ge1} 
  \sum_{\substack{n_1+n_2=n\\n_1,n_2\ge1}} \frac{(n-1)!}{ n_1! n_2!} \ n_1^{n_1-1} n_2^{n_2-1}   \ \frac{z^{n}}{ n!} \\
     &=  
  t^2 
    \sum_{n\ge1} 
  \sum_{\substack{n_1+n_2=n\\n_1,n_2\ge1}} 
  \frac{1}{ n} {n \choose \ n_1}\ n_1^{n_1-1} n_2^{n_2-1}   \ \frac{z^{n}}{ n!}
\end{align*}
Finally, Lemma \ref{le:count_planted_fixed_type} tells us that the coefficient of $\frac{z^{n}}{ n!}$ is the number  of rooted Cayley trees with two records.

 \end{ex}

The natural generalization of these examples gives the generating function for  trees with a fixed number of records (Proposition \ref{le:trees_k_records}). 

\begin{p} 
\label{le:trees_k_records}
The generating function for   trees with $j$ records can be expressed in terms of the Cayley tree function as
\begin{align*}
\int_0^z  \frac{1}{s}
  (t \mathcal{T}(s))^j  \ ds
  \end{align*}
\end{p}
\noindent
Since the number of records is fixed, the role of the variable $t$ is redundant in this Lemma. However, it will be of use in what follows (Proposition \ref{le:planted_seq}).

 \begin{proof} We start just like in the previous examples:
 \allowdisplaybreaks
\begin{align*}
\int_0^z  \frac{1}{s}
  (t \mathcal{T}(s))^j  \ ds
  &=
\int_0^z  \frac{1}{s}
  \left(t \sum_{n\ge1} n^{n-1} \frac{s^n}{n!}\right)^j  \ ds \\
  &= \int_0^z  \frac{1}{s}
  \sum_{(n_1, n_2, 
   \ldots, n_j) \vDash n} \frac{n_1^{n_1-1}}{n_1!}
   \frac{n_2^{n_2-1}}{n_2!} \ldots  
   \frac{n_j^{n_j-1}}{n_j!} t^j s^n\ ds \\
   &= 
  \sum_{(n_1, n_2, 
   \ldots, n_j) \vDash n} 
  n!\,
  \frac{n_1^{n_1-2}}{n_1!}
   \frac{n_2^{n_2-2}}{n_2!} \ldots  
   \frac{n_j^{n_j-2}}{n_j!} t^j \int_0^z  \frac{s^{n-1}}{n!}\ ds \\
    &= 
  \sum_{(n_1, n_2, 
   \ldots, n_j) \vDash n} \frac{(n-1)!}{n_1!\ n_2! \ldots n_j!} \ 
  n_1^{n_1-1} n_2^{n_2-1}  \ldots  
    n_j^{n_j-1}  t^j \ \frac{z^{n}}{n!} 
\end{align*}
where ${(n_1, n_2, 
   \ldots, n_j) \vDash n}$ denotes a composition of $n$ with $j$ parts, all positive.
   Finally,  Lemma \ref{le:count_planted_fixed_type},  (Eq.~\ref{eq_no_simpl}), implies that
the coefficient of $  \frac{z^{n}}{n!}t^j  $ in this last expression counts the number of rooted
  trees with bonsai sequence of type $t.$

\end{proof}

From this result, we derive a differential equation for the tree record function in terms of the Cayley tree function (Proposition \ref{le:planted_seq}).

\begin{p}
\label{le:planted_seq}
The  tree record generating function satisfies the differential equation
\begin{align*}
{\mathcal{R}_\bullet(z,t)} 
&= \int_0^z  \frac{1}{s}
\frac{t \, \mathcal{T}(s)   }{1-  t \, \mathcal{T}(s)} \ ds .
\end{align*}
\end{p}

\begin{proof} 
The expansion
\begin{align*}
\frac{t \mathcal{T}(s)   }{1-  t \mathcal{T}(s)} = \sum_{j \ge 1} \left( t\,  \mathcal{T}(z) \right)^j
\end{align*}
together with Proposition \ref{le:trees_k_records} imply this result. 
\end{proof} 

Combining Lemma \ref{le:Exponential_formula_record_gf}
and Proposition \ref{le:planted_seq}, we find a functional equation for the forest record function involving the Cayley tree function.

\begin{p}
\label{cor:R2B}
The forest  record generating function satisfies the differential equation:
\begin{align}
\label{eq:gf_tree2forest}
 \GFR(z,t)  = \exp \big(\int_0^z  \frac{1}{s}
\frac{t \, \mathcal{T}(s)   }{1-  t \, \mathcal{T}(s)} \ ds \big)
\end{align}

\end{p}

Finally, taking logarithms on both sides of Eq.~(\ref{eq:gf_tree2forest}), and then taking the derivative with respect to $z$, we obtain an elegant  functional equation relating records of rooted forests, and the Cayley tree function.

\begin{p}
The forest record generating function satisfies
\begin{align}
\label{eq:recordPDE}
\frac{\partial \log \GFR(z,t)}{\partial z}
& =  \frac{1}{z}
\frac{t \, \mathcal{T}(z)   }{1- t \, \mathcal{T}(z)}
\end{align}
\end{p}

Even if this result follows from the previous arguments, we give it a direct combinatorial proof. 
Note that this serves as an alternative proof of Proposition \ref{le:trees_k_records}.

\begin{proof} 
Since 
$\partial_z \log \GFR(z,t)
=
 \frac{\partial_z\GFR(z, t)}{\GFR(z,t)}
$, Eq.~(\ref{eq:recordPDE}) can be rewritten as
\[
z \partial_z\GFR(z, t) 
= 
t \, \mathcal{T}(z) \ \GFR(z, t)
+
z \partial_z\GFR(z, t)\ t \, \mathcal{T}(z)
\]
where we denote by
$\partial_z\GFR(z, t)$ the derivative of $\GFR(z,t)$ with respect to  $z$.

We claim that the generating functions appearing on both sides of this equation count  trees rooted at $\circ$ in which one of the nodes (different from $\circ$) has been selected. 
The left hand side is immediate, since
\begin{align*}
z \partial_z\GFR(z, t)  =
\sum_{n,k\ge 0} n R(n,k) \frac{z^{n}}{n!} t^k 
\end{align*}
and any   tree counted in $R(n,k)$ has
 $n$ nodes different from $\circ$. Distinguishing one of them accounts for the factor of $n$.
 
The count performed on the right-hand side is more subtle, as it first divides   trees according to whether its node with the largest label is attached to $\circ$ or not. We consider these two situations separately.

 On the one hand, the summand \emph{$t \, \mathcal{T}(z) \, \GFR(z, t)$} gives the generating function for rooted   trees where the largest node is attached to $\circ$. 
To see this, assume that $T$ is a tree rooted at $\circ$ and counted by $R(n,k)$,  and that its largest node  is attached to $\circ$. If we delete from $T$ the edge with endpoints $n$ and $\circ$, we split $T$, bijectively, into a   tree counted by $R(m_1, k)$ and a bonsai tree counted by $m_2^{m_2-2}$, where $m_1+m_2=n$. The product of the weights of these two subtrees is
\begin{align*}
\label{eq:reshuffling}
&\frac{z^{m_1}t^{k-1}}{m_1!} 
\frac{z^{m_2}t}{(m_2-1)!} 
= n \ {n-1\choose m_2-1} \ \frac{z^{n}t^{k}}{n!}.
\end{align*}
To conclude the argument, observe that the binomial coefficient appearing in the equation counts the number of ways of choosing a set of $m_2-1$
labels to relabel the nodes of the bonsai tree (the largest node will always be $m_1+m_2=n$, we only need to choose the remaining nodes). The   tree is labelled with the remaining entries. Moreover, in both cases, we keep the relative order of the nodes during the process. Finally, the factor of $n$ that appears on the right-hand side of the previous equation has the effect of distinguishing a node. 

On the other hand, the summand \emph{$z  \, \mathcal{T}(z)\GFR'(z, t)$} gives the generating function for rooted   trees where the largest node is attached to a node different from $\circ.$
 To see this, observe that $T$ is a tree counted in $R(n,k)$ but its largest node $n$ is not attached to $\circ$, we reason as follows. Since $n$ has to be attached to some other node, it has the effect of distinguishing a node. We  break $T$ into an ordered pair consisting of a   tree with a distinguished node, and a bonsai tree. 
The previous displayed equation allows us to conclude that the product of the weights of these two trees has the effect of distinguishing a node of $T.$ Therefore, these trees have as generating function $z  \, \mathcal{T}(z)\GFR'(z, t)$.

Finally, since these two situations are mutually exclusive, we obtain the desired result.
\end{proof}

\begin{p}    
\label{eq:forest_to_cayley_tree_functions}
The tree and forest record functions obey the following formulas: 
\begin{enumerate}[label=(\alph*)]
\item For fixed $k$, the generating function for the numbers $R_\bullet(n, k)$ is given by 
\[
[t^k]\mathcal R_\bullet(z, t) = \frac{\mathcal{T}(z)^k}{k} - \frac{\mathcal{T}(z)^{k+1}}{k+1}. 
\]
\item The tree record function is given by
\[
\mathcal R_\bullet(z,t) = \mathcal T(z) - \frac{t-1}{t}\log(1-t\mathcal T(z)).
\]
\item The forest record function is given by
\[
\mathcal R(z, t) = \frac{1}{z} \mathcal{T}(z)(1-t\mathcal T(z))^{-\frac{t-1}{t}}.
\]
\end{enumerate}
\end{p}
\begin{proof}
We derive (a) by directly solving the integral of Proposition \ref{le:trees_k_records},
\[
[t^k]\mathcal R_\bullet(z, t) = \int_{0}^z \frac{1}{s} \mathcal T(s)^kds.
\]
We make the change of variables $w = \mathcal T(s)$. From the functional equation $\mathcal T(s) = s \exp(\mathcal T(s))$ we have 
\[
s = w e^{-w}, \qquad ds = (e^{-w}-we^{-w})dw,
\]
and hence 
\begin{align*}
   [t^k]\mathcal R_\bullet(z, t) &=  \int_{0}^{\mathcal T(z)} \frac{1}{we^{-w}} w^k(e^{-w}-we^{-w})dw \\ &= \int_0^{\mathcal T(z)} (w^{k-1 } - w^k)dw =   \frac{T(z)^k}{k} -  \frac{T(z)^{k+1}}{k+1}.
\end{align*}
This change of variables also does the trick for (b). The integral of Proposition \ref{le:planted_seq} becomes 
\[
\begin{split}
\mathcal R_\bullet(z, t) &=\int_0^{\mathcal T(z)}\frac{1}{we^{-w}}\frac{tw}{1-tw}(e^{-w}-we^{-w})dw \\&= \int_0^{\mathcal T(z)} t\frac{1-w}{1-tw}dw=
\int_0^{\mathcal T(z)}\left(1 + \frac{t-1}{1-tw}\right)dw \\
&= \mathcal T(z) - \frac{t-1}{t}\log(1-t\mathcal T(z)).
\end{split}
\]
Finally, item (c) follows combining (b) with Lemma \ref{le:Exponential_formula_record_gf}. Indeed,  
\[
\mathcal R(z, t) = \exp(\mathcal R_\bullet(z,t)) = \exp(\mathcal T(z))(1-t\mathcal T(z))^{-\frac{t-1}{t}} = \frac{1}{z}\mathcal{T}(z)(1-t\mathcal T(z))^{-\frac{t-1}{t}}
\]
where in the last equality we have used that $z\exp(\mathcal T(z))={\mathcal{T}(z)}$.
\end{proof}

%% file: 4_record_numbers.tex
\section{Record numbers}

\begin{thm}
\label{prop:R_planted}
  The number of  rooted   trees labelled with  $[n]$ and  with $k$ records obeys the equation
    \[
    \planted(n,k) = 
       k \,(n-1)\cdots
    (n-k+1) \, n^{n-k-1}.
     \]
\end{thm}
\noindent
Remark that, for $k=1$, the product appearing in this expression has to be read as an empty product and  interpreted as 1.
    
    \begin{proof}
    The generating function for rooted forests with $k$ connected components is given by  \begin{equation}
    \label{eqn: T^k as egf}
    \frac{\mathcal T(z)^k}{k!} = \sum_{n \geq k}\binom{n}{k}kn^{n-k-1}\frac{z^n}{n!}.
    \end{equation}
 This result follows directly from Cayley’s elegant 1889 formula  \cite{Cayley}  for enumerating rooted forests with a prescribed set of roots.  
Let $S \subseteq [n]$ be a $k$-element subset.  
Then, the number of rooted forests on the vertex set $[n]$ having $k$ trees whose roots are precisely the elements of $S$ is given by
$
k\,n^{\,n - k - 1}.
$

    By Proposition \ref{le:trees_k_records} we have 
    \[\begin{split}
    \sum_{n\geq k} R_\bullet(n, k)\frac{z^n}{n!} &= \int_0^z \frac{1}{s} \mathcal{T}(s)^kds\\ &= \int_0^z \left(\sum_{n\geq k}k! \binom{n}{k}kn^{n-k-1}\frac {s^{n-1}}{n!}\right)ds \\ &= \sum_{n\geq k}k!\binom{n}{k}kn^{n-k-2}\frac{z^n}{n!}.
    \end{split}\]
    Extracting coefficients from both sides of this equation proves the claim.
\end{proof}

Note that combining the elegant formula for the tree record numbers presented in Theorem \ref{prop:R_planted}, 
with the recurrence for them presented in Proposition
\ref{le:recurrence_trees}, we obtain a nontrivial combinatorial identity.
 
\begin{cor}  For any positive integer $n$,
   \[
   n^{n-2} = \sum_{k=1}^n
       k \, \frac{n!}{(n-k)!} \, n^{n-k}.
     \]
\end{cor}

\begin{cor}
\label{cor:record_fix_n}
For fixed $n$, we have 
    \[
    \sum_{k=1}^n \planted(n,k) \frac{t^k}{k!} =  \frac{t}{n} (n+t)^{n-1}.
    \]
    \end{cor}
\begin{proof} This results can be obtained by direct calculation.
    Newton's binomial expansion yields
    \begin{align*}
    \frac{t}{n} (n+t)^{n-1}  &= \sum_{j=0}^{n-1} \binom{n-1}{j} n^{n - j - 2} t^{j+1} =
    \sum_{k=1}^n \binom{n-1}{k-1}n^{n-k-1} t^{k}\\
    &=
    \sum_{k=1}^n k(n-1)\cdots(n-k+1)n^{n-k-1} \frac{t^k}{k!} 
    = \sum_{k=1}^n R_\bullet(n,k)\frac{t^k}{k!}. \qedhere
    \end{align*}
\end{proof}

It is remarkable that the generating equation that emerges is exponential, given that previously we were using an ordinary generating function with respect to the variable $t$. 
\\

\begin{thm}
\label{thm:R(n,k) Stirling}
The number of rooted forests labelled with $[n]$ and with $k$ records obeys the equation
\[
R(n, k) = \sum_{m = k+1}^{n+1} (-1)^{m+k-1} \binom{n}{m-1}m(n+1)^{n-m}c(m, m-k).
\]    
where we denote by $c(\cdot,\cdot)$ 
the unsigned Stirling numbers of the first kind.
\end{thm}

 \begin{proof}
    The generalized binomial theorem allows us to write the identity of Proposition \ref{eq:forest_to_cayley_tree_functions} (c) as
     \begin{align*}
     \mathcal R(z, t) = \frac{\mathcal T(z)}{z} \sum_{m=0}^{\infty} \binom{-\frac{t-1}{t}}{m} (-t\mathcal T(z))^{m}
     &= \sum_{m=1}^{\infty}\frac{P_{m-1}(t)}{(m-1)!} \frac{\mathcal T(z)^m}{z}
     \end{align*}
     where $
     P_m(t) = (t-1)(2t-1)\cdots(mt-1).
     $
     Thus, the coefficient of $z^n$ in the forest record generating function is given by
     \[
     [z^n]\mathcal R(z, t) = \sum_{m= 1}^{\infty} \frac{P_{m-1}(t)}{(m-1)!} [z^{n+1}]\mathcal T(z)^m = \frac{1}{n!}\sum_{m=1}^{n+1} 
    \binom{n}{m-1}m (n+1)^{n-m}P_{m-1}(t),
     \]
     where the last equality follows from Equation \eqref{eqn: T^k as egf}.
     Therefore,
     \begin{equation}
     \label{eqn: R(n,k) in terms of P_m}
     R(n,k) =n![t^kz^n] \mathcal R(z, t) = \sum_{m=k+1}^{n+1} \binom{n}{m-1}m(n+1)^{n-m} [t^k]P_{m-1}(t). \end{equation}
    To conclude the proof, note that 
    \[
    P_{m-1}(t) = \prod_{i=1}^{m-1} (it - 1) = (-1)^{m-1}t^m \prod_{i=0}^{m-1} (\frac{1}{t} - i) = (-1)^{m-1}t^m \left(\frac{1}{t}\right)_{m}, 
    \]
    where $(\cdot)_m$ denotes the falling factorial. 
    Recall that $(x)_m = \sum_{k=0}^m (-1)^{m-k}c(m,k)x^k$, and hence 
    \[
    P_{m-1}(t) = \sum_{k=1}^{m-1}(-1)^{m-1 + m-k} c(m, k)t^{m-k} = \sum_{k=1}^{m-1}(-1)^{m+k-1}c(m, m-k)t^{k}.
    \]
    In particular, the coefficient of $t^k$ in this polynomial is  \[[t^k] P_{m-1}(t) = (-1)^{m+k-1}c(m, m-k).\]
    Combining this expression with Equation \eqref{eqn: R(n,k) in terms of P_m}
    yields the desired result.
 \end{proof}
 
\begin{cor} The forest record numbers can be expressed in terms of the tree record numbers as 
    \[
    R(n, k) = \sum_{m=k+1}^{n+1} \frac{(-1)^{m+k-1}}{(m-1)!}c(m, m-k) R_\bullet(n+1, m).
    \]
\end{cor}
\begin{proof}
    Combine Theorem \ref{thm:R(n,k) Stirling} with Theorem \ref{prop:R_planted}.
\end{proof}

%% file: 5_consequences.tex
\section{Consequences}

\subsection{Asymptotics}
We compute the asymptotics for the generating function for records of rooted trees. 

\begin{cor}
\label{Cor:Lemma_k_asymptotics_records}
For $k$ fixed, $\planted(n, k)$ is asymptotic to $kn^{n-2}.$
\end{cor} 
\begin{proof}
Recall that the Stirling approximation states that
    $m! \sim \sqrt{2\pi m}(m/e)^m$. Thus,     $\planted(n, k) $ satisfies that
    \begin{align*}
    \hspace{1.5cm}
        \frac{kn^{n-k-1}(n-1)!}{(n-k)!} \sim  kn^{n-k-1} \frac{\sqrt{2\pi(n-1)}\left(\frac{n-1}{e}\right)^{n-1}}{\sqrt{2\pi(n-k)}\left(\frac{n-k}{e}\right)^{n-k}} 
        =  kn^{n-k-1} \sqrt{\frac{n-1}{n-k}}\frac{(n-1)^{n-1}}{(n-k)^{n-k}}e^{1-k}.
    \end{align*}
Moreover, for $k$ fixed and $n$ tending to infinity, we find that
\begin{align*}
 \hspace{1.5cm}
    (n-1)^{n-1} \sim 
    \frac{n^{n-1}}{e} \  \text{ and }\  
    (n-k)^{n-k} \sim 
     \frac{n^{n-k}}{e^k}.
    \end{align*}
    Finally,
    \[ 
    \sqrt{\frac{n-1}{n-k}} \sim 1,
\ \text{and so, }\
    kn^{n-k-1} \sqrt{\frac{n-1}{n-k}}\frac{(n-1)^{n-1}}{(n-k)^{n-k}}e^{1-k} \sim kn^{n-2}. \qedhere
\] 
\end{proof}

\subsection{Log-concavity} The tree record numbers are log-concave. This result is easy to obtain out of their elegant formula. On the other hand, we have verified using SageMath that the forest record numbers  remain log-concave up to \emph{\( n = 900 \)}.

\begin{cor} 
\label{cor:log_concavity_trees}
For each natural number $n$, the sequence  $\planted(n, k)$, for $k\ge 1$ is log-concave.
\end{cor}
\begin{proof} In Theorem \ref{prop:R_planted}, we showed that this sequence is the termwise product of positive integers
$k!\binom{n}{k}$ and $kn^{n-k-2}$, both of which are log-concave on $k$. \end{proof}

We are quoting a result that is folklore, but we have not been able to find a precise reference for it.  We add a proof. Two sequences of positive integers $(a_k)_k$ y $(b_k)_k$ are log-concave if and only if $a_k^2 \geq a_{k-1}a_{k+1}$ and $b_k^2 \geq b_{k-1}b_{k+1}$ for all $k$. Therefore, the ``term-wise product" verifies that for all $k$: $ (a_kb_k)^2 \geq (a_{k-1}b_{k-1})(a_{k+1}b_{k+1})
$. 

In 1968, E.-H. Lieb \cite{Lieb} proved that both families of Stirling numbers are log-concave showing that their generating functions have only real roots.  However, the polynomials obtained from the rows of the table of forest record numbers (Table \ref{fig:planted_table}) can have non-real roots. For example, from the second row we obtain \(3q + 7q^2 + 6q^3\), a polynomial that has a pair of nonreal roots.  However, this observation does not contradict the possibility that the record numbers are log-concave. 
Recall that log-concave sequences are in particular unimodal, but unimodality tends to be harder to show directly. 
\\

\subsection{Asymptotic peak location}
We determine the value of \(k\) that maximizes the tree record numbers \(R_\bullet(n,k)\) for a fixed \(n\). We then study the analogous problem for the forest record numbers.

\begin{table}[h!]
\centering
     \resizebox{0.9\textwidth}{!}{\input{tikz.peak_location}
     \input{tikz.peak_location_squared}}
\caption{On the left, the index maximizing $R(n,k)$ for $n\le 900$. On the right, the same index, but squared. The least squares approximation of the data on the right has slope 1.2158 and intercept 35.7936.}
\label{Fi:peak_fits}
 \end{table}
 
\begin{cor}
\label{co:location_peaks_trees}
    For fixed $n$, the greatest index $k$ maximizing $R_\bullet(n, k)$ is given by 
    \[
    k = \left\lfloor \frac{1 + \sqrt{1 + 4n}}{2} \right\rfloor.
    \]
\end{cor}
\begin{proof}
    Case $n=1$ trivially holds, so let $n \geq 2$. By Corollary \ref{cor:log_concavity_trees}, the sequence $R_\bullet(n, k)$ is log-concave, and thus unimodal. This implies that the number $R_\bullet(n, k)$ is maximised just before the first index $k$ for which the quotient $\frac{R_\bullet(n, k)}{R_\bullet(n, k-1)}$ is strictly lower than 1. But 
    \[
    \frac{R_\bullet(n, k)}{R_\bullet(n, k-1)} = \frac{(n-1)\cdots(n-k+1)kn^{n-k-1}}{(n-1)\cdots(n-k + 2) (k-1)n^{n-k}} = \frac{(n-k+1)k}{(k-1)n}.
    \]
    Direct computation shows that this quotient is greater than or equal to 1 when $1 < k \leq \frac{1}{2} (1 + \sqrt{1+4n})$, and lower than 1 for any greater values of $k$. Thus, the index $k$ maximizing $R_\bullet(n, k)$ is attained at the greatest integer $k \leq \frac{1}{2}(1+\sqrt{1+4n})$.
\end{proof}

The value of $k$ maximizing $R(n,k)$,  the forest record numbers, appears to grow on the order of $\sqrt{n}$. See the graphs depicted in Table \ref{Fi:peak_fits}.
It remains to determine the precise asymptotic position of the peak locations.
 Compare this  with the analogous situation for the Stirling numbers, the index \( k \) maximizing \( c(n, k) \) satisfies
that $
  k \sim \log n \quad \text{as } n \to \infty$
  matching the expected number of cycles in a random permutation.

%% file: tikz.peak_location.tex
\begin{tikzpicture}
  \begin{axis}[
    width=10cm,
    height=8cm,
    xlabel={$n$},
    ylabel={$k$},
    grid=both,
    grid style={line width=.1pt, draw=gray!10},
    major grid style={line width=.2pt,draw=gray!35!white},
    legend style={at={(0.54,0.36)}, anchor=north west, font=\small},
    tick label style={font=\small},
    label style={font=\small},
    title style={font=\small},
    ymin =0,
    ymax=35,
    xmin=0,
    xmax=900
  ]

    \addplot graphics [xmin=1, xmax=900, ymin=-5, ymax=35]
      {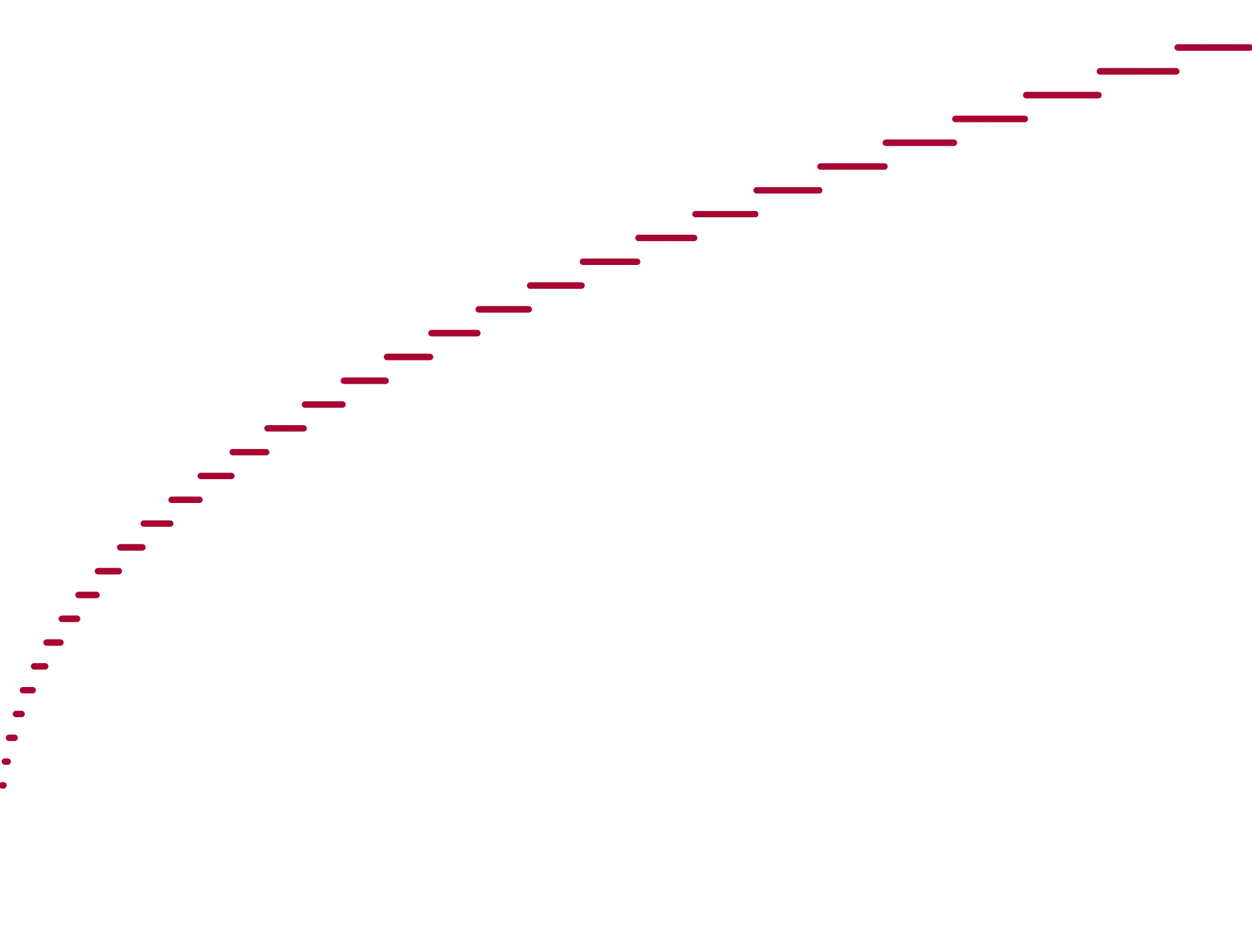}; 



  \end{axis}
\end{tikzpicture}

%% file: tikz.peak_location_squared.tex
\begin{tikzpicture}
  \begin{axis}[
    width=10cm,
    height=8cm,
    xlabel={$n$},
    ylabel={$k$},
    grid=both,
    grid style={line width=.1pt, draw=gray!10},
    major grid style={line width=.2pt,draw=gray!35!white},
    legend style={at={(0.54,0.36)}, anchor=north west, font=\small},
    tick label style={font=\small},
    label style={font=\small},
    title style={font=\small},
    ymin =0,
    ymax=1200,
    xmin=0,
    xmax=900
  ]

    \addplot graphics [xmin=1, xmax=900, ymin=0, ymax=1200]
      {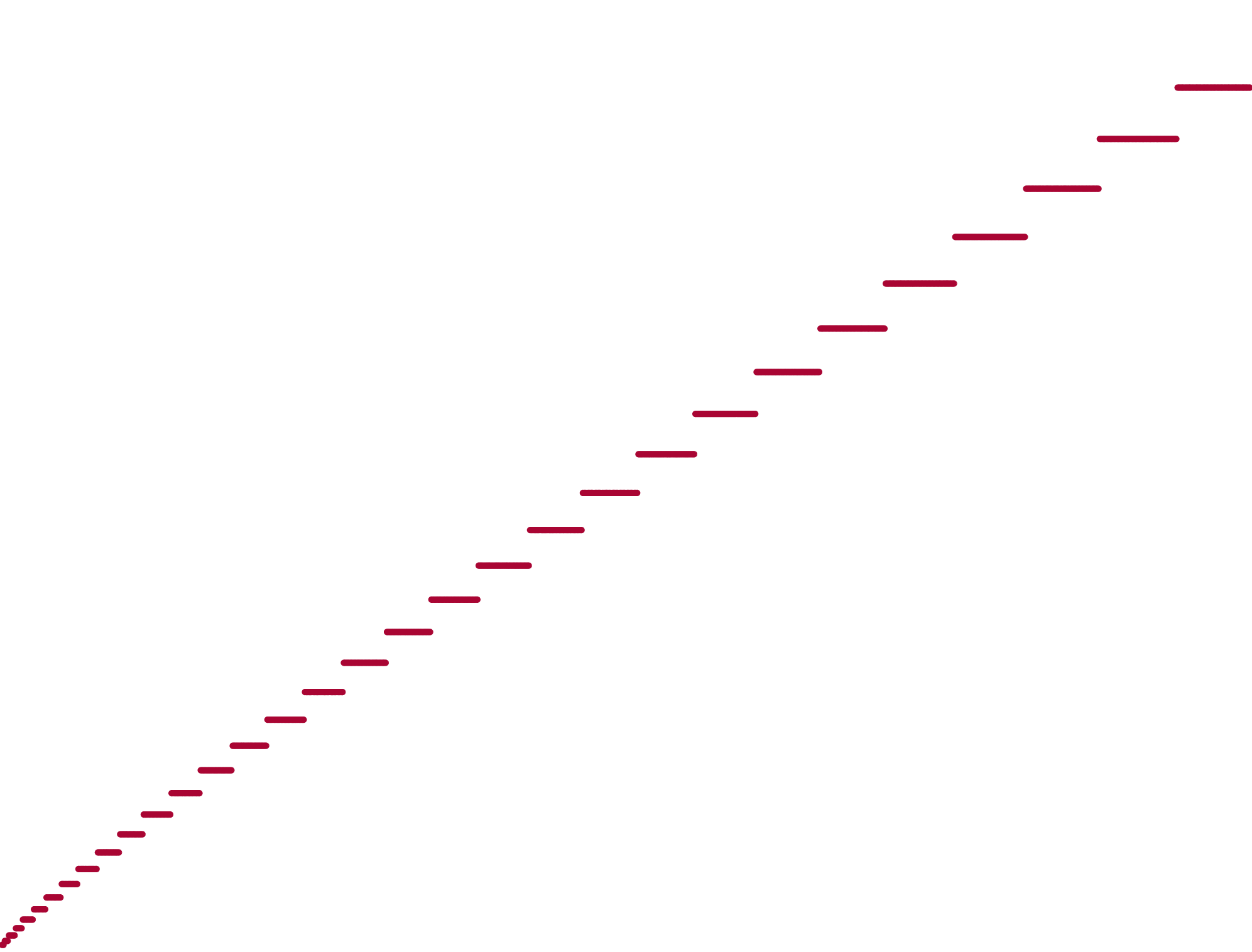}; 



  \end{axis}
\end{tikzpicture}

%% file: 6_final_comments.tex
\section{Final Comments}

\begin{enumerate}[wide, labelindent=0pt]

\item 
Most constructions presented in this work have been implemented in Sage \cite{sage}, and the corresponding source code can be found in \cite{RecordsSage}. 

\item  I. Kortchemski introduced the record-sum polynomials  of  the symmetric group. He defined $\mathcal{C}(n,k)$ as the number of permutations in $\SS_n$ for which the sum of the sum of the positions of
their records is equal to $k$. Then, Kortchemski showed (Proposition 2.6, \cite{Kortchemski}) that its generating function is
\[
\sum_{\sigma \in \SS_n} q^{\text{srec}(\sigma)} = q (q^2 + 1)(q^3 + 2) \cdots (q^n + n - 1) = \sum_{k=1}^{\frac{n(n+1)}{2}} \mathcal{C}(n,k) q^k.
\]
We are interested in understanding the natural analogue of Kortchemski's statistic for trees, where instead of considering the positions of the records, we use their depth. (Recall that the depth  of a node $v$ is the number of edges present in path from the root to $v$.)

\[
P_n(q) = \sum_{T \text{   tree}} q^{\srec(T)} 
\]

\input{table.srec}

The coefficients of the firsts record-sum polynomials for rooted   trees are tabulated in Table \ref{fig:srecord_table}.  It appears that the second column of this table is \url{https://oeis.org/A083483}: Number of (unrooted) forests with two connected components and $n$ vertices. 
On the other hand, the rows of this table are not unimodal, which in turns implies that they are not log-concave.\\

\item 
Numbers counted by the formula for records of rooted trees presented in Theorem \ref{prop:R_planted} have already appeared in the literature of queues and congestion theory. Building on the work of Borel \cite{borel1942emploi}, Tanner \cite{Tanner} established that $\frac{\planted(n, k)}{(n-1)!} e^{-\rho n} \rho^{n-k}$ represents (in the words of Haight \cite{Haight}) ``{\it the probability that exactly $k$ members of a queue will be served before the queue first becomes empty, beginning with $n$ members and with traffic intensity $\rho$, assuming Poisson input and regular discharge}". It would be interesting to further understand the interplay between records in rooted trees and the Borel-Tanner distribution. \\

 \item If $\mathcal{A}$ and $\mathcal{B}$ be combinatorial classes with exponential generating functions $A(z)$ and $B(z)$.
Then the product $A(z)B(z)$ is precisely the exponential generating function of the labelled product $\mathcal{A} \times \mathcal{B}$. If we add the condition that the largest label is constrained to lie in the $\mathcal{A}$ component, we define a new product, the \emph{boxed product}. The generating function of the resulting class is:
\[
\int_0^{z} \left(\partial{A(t)}\right) B(t) \ dt
\]
More information on the boxed product can be found in \cite{Flajolet}, where the boxed product is used to obtain an elegant derivation for  the generating function for record in permutations (Example II.16). However,  while the situation for Cayley trees is closely related, is it definitely different because there is not an unique way of reconstructing a Cayley tree from its Bonsai sequence.\\
 
  \item An interesting question, raised by Cyril Bandelier, was whether the expressions for records of trees and forests remain valid when the Cayley tree function is replaced by other functions (for example, the generating function for binary trees). This does not appear to be the case. The reason seems to be that the other families of trees we have tried are not compatible with the record decompositions. In particular, a tree obtained from a bonsai sequence of a binary tree is not necessarily binary.

 \end{enumerate}

%% file: table.srec.tex

\begin{table}[h!]
\renewcommand{\arraystretch}{1.3} 
\rowcolors{2}{gray!20}{white} 

    \begin{tabular}{c|ccccccccccccccccccccc}
        \toprule
$n \backslash k$ & $1$ & $2$ & $3$ & $4$ & $5$ & $6$ & $7$ & $8$ & $9$ & $10$ & $11$  & $12$  & $13$  & $14$  & $15$
        \\
        \midrule
$1$ & $1$   & $0$   & $0$   & $0$  & $0$  & $0$  & $0$  & $0$  & $0$  & $0$  & $0$  & $0$  & $0$  & $0$  & $0$ \\
$2$ & $1$   & $1$   & $1$   & $0$  & $0$  & $0$  & $0$  & $0$  & $0$  & $0$  & $0$  & $0$  & $0$  & $0$  & $0$ \\
$3$ & $3$   & $3$   & $4$   & $4$   & $1$   & $1$   & $0$  & $0$  & $0$  & $0$  & $0$  & $0$  & $0$  & $0$  & $0$ \\
$4$ & $16$  & $15$  & $21$  & $26$  & $18$  & $15$  & $8$   & $4$  & $1$  & $1$   & $0$  & $0$  & $0$  & $0$  & $0$\\
$5$ & $125$ & $110$ & $155$ & $203$ & $177$ & $173$ & $135$ & $94$ & $59$ & $31$ & $20$ & $8$ & $4$ & $1$ & $1$ \\

        \bottomrule
    \end{tabular}
        \caption{The number of trees with $n$ vertices such that the sum of the depth of its records is $k$.} 
        \label{fig:srecord_table}
\end{table}
\rowcolors{2}{}{} 